\documentclass[12pt]{amsart}
\usepackage{amsfonts}
\usepackage{amssymb}
\usepackage{amscd}
\usepackage{color}

\newtheorem{theorem}{Theorem}[section]
\newtheorem{lemma}[theorem]{Lemma}

\newtheorem{proposition}[theorem]{Proposition}

\theoremstyle{definition}
\newtheorem{definition}[theorem]{Definition}

\newtheorem{example}[theorem]{Example}
\theoremstyle{remark}
\newtheorem{remark}[theorem]{Remark}

\newtheorem{note}[theorem]{Note}

\begin{document}
\title{Turing degrees of isomorphism types of algebraic objects}
\author{Valentina Harizanov}
\address{Department of Mathematics\\
George Washington University\\
Washington, DC 20052}
\email{harizanv@gwu.edu}
\author{Wesley Calvert}
\address{Department of Mathematics \& Statistics\\
Murray State University\\
Murray, KY 42071}
\email{wcalvert@sigmaxi.org}
\author{Alexandra Shlapentokh }
\address{Department of Mathematics\\
East Carolina University\\
Greenville, NC 27858}
\email{shlapentokha@mail.ecu.edu}

\begin{abstract}
The Turing degree spectrum of a countable structure $\mathcal{A}$ is the set
of all Turing degrees of isomorphic copies of $\mathcal{A}$. The Turing
degree of the isomorphism type of $\mathcal{A}$, if it exists, is the least
Turing degree in its degree spectrum. We show there are countable fields,
rings, and torsion-free abelian groups of arbitrary rank, whose isomorphism
types have arbitrary Turing degrees. We also show that there are structures
in each of these classes whose isomorphism types do not have Turing degrees.
\end{abstract}

\thanks{Calvert acknowledges partial support of the NSF grants DMS-0139626
and DMS-0353748. Harizanov gratefully acknowledges support of the Columbian
Research Fellowship of the George Washington University. Shlapentokh
acknowledges partial support of the NSF grant DMS-0354907. }
\maketitle

\section{Introduction}

One of the main goals of computable algebra is to understand how algebraic
properties of structures interact with their computability-theoretic
properties. While in algebra and model theory isomorphic structures are
often identified, in computable model theory they can have very different
algorithmic properties. Here, we study Turing degrees of algebraic
structures from some well-known classes. We consider only countable
structures for computable (usually finite) languages. The universe $A$ of an
infinite countable structure $\mathcal{A}$ can be identified with the set $%
\omega $ of all natural numbers. Furthermore, we often use the same symbol
for the structure and its universe. (For the definition of a language and a
structure see p. 8 of \cite{Marker}, and for a definition of a computable
language see p. 509 of \cite{Millar}.)

Let $\{\mathcal{A}_{j},j\in \omega \}$ be a sequence of structures contained
in a structure $\mathcal{B}$. Then by $\prod_{j\in \omega }\mathcal{A}_{j}$
we mean the smallest substructure of $\mathcal{B}$ containing $\mathcal{A}%
_{j}$ for all $j$. More specifically, we will be looking at products of
number fields and function fields, and at products of rings contained in a
number field. In the case of fields, we fix an algebraic closure of ${\mathbb Q}$, a
rational function field over a finite field of characteristic $p>0$ or over $%
\mathbb{Q}$ as required, and we can set $\mathcal{B}$ to be this algebraic
closure. In the case of subrings of a number field, the number field itself
is a natural choice for $\mathcal{B}$.

When measuring complexity of structures, we identify them with their atomic
diagrams. The atomic diagram of a structure $\mathcal{A}$ is the set of all
quantifier-free sentences in the language of $\mathcal{A}$ expanded by
adding a constant symbol for every $a\in \mathnormal{A}$, which are true in $%
\mathcal{A}$. The \emph{Turing degree} of $\mathcal{A}$, $deg(\mathcal{A})$,
is the Turing degree of the atomic diagram of $\mathcal{A}$. Hence, $%
\mathcal{A}$ is \emph{computable }($\emph{recursive}$) iff $deg(\mathcal{A})=%
\mathbf{0}$. (Some authors call a structure computable if it is only
isomorphic to a computable one.) We also say that a procedure is \emph{%
computable} (effective), \emph{relative to} $\mathcal{B}$, also \emph{in} $%
\mathcal{B}$, if it is computable relative to the atomic diagram of $%
\mathcal{B}$\textit{. }

We use $\leq _{T}$ for Turing reducibility and $\equiv _{T}$ for Turing
equivalence of sets. We often abbreviate $deg(\mathcal{A})\leq deg(X)$ by $%
\mathcal{A}\leq _{T}X$. (Detailed information about Turing degrees and their
structure can be found in \cite{Rogers} and \cite{S}.)

The \emph{Turing degree spectrum} of a countable structure $\mathcal{A}$ is
\begin{equation*}
DgSp(\mathcal{A})=\{\deg (\mathcal{B}):\mathcal{B}\cong \mathcal{A}\}\text{.}
\end{equation*}%
A countable structure $\mathcal{A}$ is \emph{automorphically trivial} if
there is a finite subset $X$ of the domain $A$ such that every permutation
of $A$, whose restriction on $X$ is the identity, is an automorphism of $%
\mathcal{A}$. If a structure $\mathcal{A}$ is automorphically trivial, then
all isomorphic copies of $\mathcal{A}$ have the same Turing degree. It was
shown in \cite{HM} that if the language is finite, then that degree must be $%
\mathbf{0}$. On the other hand, Knight \cite{K} proved that for an
automorphically nontrivial structure $\mathcal{A}$, we have that $DgSp(%
\mathcal{A})$ is closed upwards, that is, if $\mathbf{b}\in DgSp(\mathcal{A}%
) $ and $\mathbf{d}>\mathbf{b}$, then $\mathbf{d}\in DgSp(\mathcal{A})$.
Hirschfeldt, Khoussainov, Shore, and Slinko established in \cite{HKSS} that
for every automorphically nontrivial structure $\mathcal{A}$, there is a
symmetric irreflexive graph, a partial order, a lattice, a ring, an integral
domain of arbitrary characteristic, a commutative semigroup, or a $2$-step
nilpotent group whose degree spectrum coincides with $DgSp(\mathcal{A})$.

The Turing degree of a structure is not invariant under isomorphisms. Thus,
Jockusch and Richter introduced the following complexity measures of the
isomorphism type of a structure.

\begin{definition}
$($i$)$ The \emph{Turing degree of the isomorphism type} of $\mathcal{A}$,
if it exists, is the least Turing degree in $DgSp(\mathcal{A})$.

$($ii$)$ Let $\alpha $ be a computable ordinal. The $\alpha $\emph{th jump
degree} of a structure $\mathcal{A}$ is, if it exists, the least Turing
degree in%
\begin{equation*}
\{deg(\mathcal{B})^{(\alpha )}:\mathcal{B}\cong \mathcal{A}\}\text{.}
\end{equation*}
\end{definition}

\noindent Obviously, the notion of the $0$th jump degree of $\mathcal{A}$
coincides with the notion of the degree of the isomorphism type of $\mathcal{%
A}$. (A general discussion of the jump operator can be found in 13.1 of \cite%
{Rogers} and Chapter III\ of \cite{S}.)

In \cite{R} Richter proved that if $\mathcal{A}$ is a structure without a
computable copy and satisfies the effective extendability condition
explained below, then the isomorphism type of $\mathcal{A}$ has no degree.
Richter's result uses a minimal pair construction. Distinct nonzero Turing
degrees $\mathbf{a}$ and $\mathbf{b}$ form a \emph{minimal pair} if
\begin{equation*}
(\mathbf{c}\leq \mathbf{a}\text{,\ }\mathbf{c}\leq \mathbf{b)}\Rightarrow
\mathbf{c}=\mathbf{0}\text{.}
\end{equation*}
(See \cite{S} for the minimal pair method.) A structure $\mathcal{A}$
satisfies the effective extendability condition if for every finite
structure $\mathcal{M}$ isomorphic to a substructure of $\mathcal{A}$, and
every embedding $\sigma $ of $\mathcal{M}$ into $\mathcal{A}$, there is an
algorithm that determines whether a given finite structure $\mathcal{N}$
extending $\mathcal{M}$ can be embedded into $\mathcal{A}$ by an embedding
extending $\sigma $. In \cite{R} Richter also showed that every linear order
and every tree, as a partially ordered set, satisfy the effective
extendability condition. Recently, Khisamiev in \cite{Kh} showed that every
abelian $p$-group, where $p$ is a prime number, satisfies the effective
extendability condition. Hence the isomorphism type of a countable linear
order, a tree, or an abelian $p$-group, which is not isomorphic to a
computable one, does not have a degree.

Recently, Csima in \cite{Csima} proved that if $\mathcal{A}$ is a prime
model of a complete decidable theory with no computable prime model, then
the isomorphism type of $\mathcal{A}$ does not have a Turing degree, while
for every $n\geq 1$, the structure $\mathcal{A}$ has the $n$th jump degree $%
\mathbf{0}^{(n)}$. (See \cite{H} for more information on computability of
prime models.)

If $\mathcal{A}$ is a nonstandard model of Peano arithmetic, then the
isomorphism type of $\mathcal{A}$ has no degree. On the other hand, Knight
showed in \cite{K} that for any Turing degree $\mathbf{d}$, there is a
nonstandard model of Peano Arithmetic with first jump degree $\mathbf{d}%
^{\prime }$. Knight also established that the only possible first jump
degree for a linear order is $\mathbf{0}^{\prime }$.

Ash, Jockusch, and Knight in \cite{AJK}, and Downey and Knight in \cite{DK}
proved that for every computable ordinal $\alpha \geq 1$, and every Turing
degree $\mathbf{d}$ such that $\mathbf{d}\geq \mathbf{0}^{(\alpha )}$, there
is a linear order $\mathcal{L}$ whose $\alpha $th jump degree is $\mathbf{d}$%
, and such that $\mathcal{L}$ does not have $\beta $th jump degree for any $%
\beta <\alpha $. Jockusch and Soare proved in \cite{JS1} that for a Turing
degree $\mathbf{d}$ and $n\in \omega $, if a Boolean algebra $\mathcal{B}$
has $n$th jump degree $\mathbf{d}$, then $\mathbf{d=0}^{(n)}$. They also
showed that if $\mathbf{d}\geq \mathbf{0}^{(\omega )}$, then there is a
Boolean algebra with $\omega $th jump degree $\mathbf{d}$. Oates in \cite{O}
proved that for every computable ordinal $\alpha \geq 1$, and every Turing
degree $\mathbf{d}$ such that $\mathbf{d}>\mathbf{0}^{(\alpha )}$, there is
an abelian group $\mathcal{G}$ whose $\alpha $th jump degree is $\mathbf{d}$%
, and $\mathcal{G}$ does not have $\beta $th jump degree for any $\beta
<\alpha $.

For additional background information on computability (recursion) theory,
see \cite{Rogers} or \cite{S}. For computable model theory, see \cite{D},
\cite{H}, \cite{H1}, and \cite{Rabin}. In the sections that follow we will
use some facts from algebra and number theory. The relevant algebraic
material can be found in \cite{A}, \cite{C}, \cite{Fried}, \cite{Fuchs},
\cite{Januz}, \cite{L}, \cite{Riben}, or \cite{Sh}.

\section{Turing degrees of the isomorphism types of structures}

We would like to further investigate Turing degrees of the isomorphism types
of abelian groups, rings, and fields. The following result is a modification
of Richter's Theorem 2.1 on p. 725 in \cite{R}. Richter used her theorem to
show that for every Turing degree $\mathbf{d}$, there is an abelian torsion
group whose isomorphism type has degree $\mathbf{d}$.

\begin{theorem}
\label{thm:1} Let $\mathcal{C}$ be a class of countable structures in a
finite language $L$, closed under isomorphisms. Assume that there is a
computable sequence $\{\mathcal{A}_{i},i\in \omega \}$ of computable $($%
possibly infinite$)$ structures in $\mathcal{C}$ satisfying the following
conditions.

\begin{itemize}
\item There exists a finitely generated structure $\mathcal{A}\in \mathcal{C}
$ such that for all $i\in \omega $, we have that $\mathcal{A}\subset
\mathcal{A}_{i}$.

\item For any $X\subseteq \omega $, there is a structure $\mathcal{A}_{X}$
in $\mathcal{C}$ such that $\mathcal{A}\subset \mathcal{A}_{X}$ and
\begin{equation*}
\mathcal{A}_{X}\leq _{T}X\text{,}
\end{equation*}%
and for every $i\in \omega $, there exists an embedding $\sigma $ such that
\begin{equation*}
\sigma :\mathcal{A}_{i}\hookrightarrow \mathcal{A}_{X}\text{, \ }\sigma _{|%
\mathcal{A}}=id\text{,}
\end{equation*}%
iff $i\in X$.

\item Suppose a structure $\mathcal{B}$ is such that for some $X\subseteq
\omega $, we have that $\mathcal{B\simeq A}_{X}$ under isomorphism $\tau :%
\mathcal{A}_{X}\longleftrightarrow \mathcal{B}$. Consider the set $\Lambda $
of pairs $(i,j)$ such that exactly one of the structures $\mathcal{A}_{i}$
and $\mathcal{A}_{j}$ is embeddable in $\mathcal{B}$ under an embedding $%
\sigma $ such that $(\tau ^{-1}\circ \sigma )_{|\mathcal{A}}=id$. Then there
is a procedure, computable in $\mathcal{B}$ $($that is, in the atomic
diagram of $\mathcal{B})$, which decides for every $(i,j)\in \Lambda $,
which of $\mathcal{A}_{i}$ and $\mathcal{A}_{j}$ embeds in $\mathcal{B}$.
\end{itemize}

Then for every Turing degree $\mathbf{d}$, there is a structure in $\mathcal{%
C}$ whose isomorphism type has degree $\mathbf{d}$.
\end{theorem}

\proof%

Let $D\subseteq \omega $ be such that $deg(D)=\mathbf{d}$. As usual, let
\begin{equation*}
D\oplus \overline{D}=_{def}\{2n:n\in D\}\cup \{2n+1:n\notin D\}\text{.}
\end{equation*}%
We will show that $\mathcal{A}_{D\oplus \overline{D}}$ is a structure in $%
\mathcal{C}$, whose isomorphism type has Turing degree $\mathbf{d}$.
Clearly, by assumption,
\begin{equation*}
\deg (\mathcal{A}_{D\oplus \overline{D}})\leq \deg (D\oplus \overline{D})=%
\mathbf{d}\text{.}
\end{equation*}

Now, let a structure $\mathcal{B}$ be such that there exists an isomorphism $%
\tau :\mathcal{A}{}_{D\oplus \overline{D}}\longleftrightarrow \mathcal{B}$.
We then have, by the definition of $D\oplus \overline{D}$ and an assumption
of the theorem, that for every $j\in \omega ,$%
\begin{equation*}
j\in D\Leftrightarrow (\exists \sigma )[\sigma :\mathcal{A}%
_{2j}\hookrightarrow \mathcal{B}\wedge (\tau ^{-1}\circ \sigma )_{|{}%
\mathcal{A}}=id]\text{,}
\end{equation*}%
or, equivalently,%
\begin{equation*}
j\notin D\Leftrightarrow (\exists \sigma )[\sigma :\mathcal{A}%
_{2j+1}\hookrightarrow \mathcal{B}\wedge (\tau ^{-1}\circ \sigma )_{|{}%
\mathcal{A}}=id]\text{.}
\end{equation*}%
Thus, by an assumption of the theorem, we conclude that $\mathbf{d\leq }\deg
(\mathcal{B})$. Hence $\deg (\mathcal{A}_{D\oplus \overline{D}})=\mathbf{d}$%
. Moreover, the degree of the isomorphism type of $\mathcal{A}_{D\oplus
\overline{D}}$ is $\mathbf{d}$.
\endproof%

\section{Fields whose isomorphism types have arbitrary Turing degrees}

In this section we show that certain sequences of computable linearly
disjoint fields satisfy conditions of Theorem \ref{thm:1}. First we state
two definitions to describe sequences we will consider.

\begin{definition}
Let $F$ be any field. Let $\{L_{i},i\in \omega \}$ be a sequence of finite
extensions of $F$ such that for any $i\in \omega $, we have that $L_{i}$ and
$\prod_{j\in \omega {}\setminus \{i\}}L_{j}$ are linearly disjoint over $F$,
or, in other words,
\begin{equation*}
\lbrack L_{i}:F]=[L:\prod_{j\in \omega {}\setminus \{i\}}L_{j}]>1,
\end{equation*}%
where $L=\prod_{j\in {}\omega }L_{j}$. Then we will call the sequence $%
\{L_{i},i\in \omega {}\}$ \textit{{totally linearly disjoint} over }$F$%
\textit{. }
\end{definition}

\begin{definition}
Let $F$ be any field. Let $\{L_{i},i\in \omega \}$ be a sequence of
algebraic extensions of $F$, and let $L=\prod_{i\in {}\omega }L_{i}$.
Suppose further that for any embedding $\sigma :L\hookrightarrow \tilde{F}$,
where $\tilde{F}$ is the algebraic closure of $F$, such that $\sigma
_{|{}F}=id$, for all $i$, we have either $\sigma (L_{i})=L_{i}$ or $\sigma
(L_{i})\not\subset $ $L$. Then we will call the sequence $\{L_{i},i\in
\omega {}\}$ \emph{stable with respect to} $F$. If $F=\mathbb{Q}$ or $F$ is
a finite field, then we will say that the sequence $\{L_{i},i\in \omega \}$
is \emph{stable}.
\end{definition}

We are now ready to prove the main theorem involving these field sequences.

\begin{theorem}
\label{thm:2}Let $K$ be any computable finitely generated field, and let $%
\tilde{K}$ be a computable algebraic closure of $K$. Let $\mathcal{P}%
=\{f_{i}(T)\in K(T)\}$ be a computable sequence of monic polynomials
irreducible over $K$. Let $\alpha _{i}$ be a root of $f_{i}$, and let $%
M_{i}=K(\alpha _{i})$. Assume further that the sequence $\{M_{i},i\in \omega
\}$ is totally linearly disjoint over $K$ and is stable with respect to $K$.
Let $\mathcal{A}=K$. Let $\mathcal{A}_{i}=M_{i}$, and for any $X\subset
\omega $, let $\mathcal{A}_{X}=M_{X}$, where
\begin{equation*}
M_{X}=\prod_{i\in X}M_{i}.
\end{equation*}%
Then the conditions of Theorem \ref{thm:1} are satisfied.
\end{theorem}

\begin{remark}
Before we proceed with the details of the proof, we would like to note that
every computable field has a computable algebraic closure, as shown in \cite%
{Rabin}.
\end{remark}

\begin{proof}
Let $X\subset \omega $.

\begin{enumerate}
\item We show that $M_{i}$ is embeddable into $M_{X}=\prod_{i\in X}M_{i}$
under any embedding $\sigma $ keeping $K$ fixed if and only if $i\in X$. Let
$\sigma :M_{i}\hookrightarrow M_{X}$ be such an embedding. Given our
assumptions on the linear disjointness of the elements of the sequence, $%
\sigma $ can be extended to $\tilde{\sigma}:M\hookrightarrow \tilde{F}$,
where $M=\prod_{j\in \omega }M_{j}$. Thus, by our assumptions, we conclude
that $\sigma (M_{i})=M_{i}$, and therefore, $M_{i}\subset M_{X}$. However,
by our assumptions again, this is possible iff $i\in X$.

\item Next, we show that $M_{X}\leq _{T}X$. Let $\alpha _{i}$ be a root of $%
f_{i}$. Given our assumptions, $M$ is a computable field, and there exists a
computable function which, given an element of $M$, will produce its
coordinates with respect to the basis%
\begin{equation*}
\Omega =\{\prod_{i\in I}\alpha _{i}^{a_{i}}:I\subset \omega \text{, \ }%
|I|<\infty \text{, \ }0\leq a_{i}<\deg (f_{i})\}\text{.}
\end{equation*}%
Thus, given $\beta \in M$ and $X$, we can determine, computably in $X$,
whether $\beta \in M_{X}$. Consequently, $M_{X}\leq _{T}X$. (Since the
reverse reducibility is obvious, we can actually show that $M_{X}\equiv
_{T}X $.)

\item Finally, let $\tau :M_{X}\longleftrightarrow \mathcal{B}$ be an
isomorphism, and let $i,j\in \omega $, $i\not=j$, be such that for exactly
one of $k=i$ and $k=j$, there exists $\sigma :M_{k}\hookrightarrow \mathcal{B%
}$ satisfying the condition that $\tau ^{-1}\circ \sigma $ is an identity on
$K$. Since $K$ is finitely generated, we can compute in $\mathcal{B}$ the $%
\tau $-images of $f_{i}$ and $f_{j}$. Next we note that there exists an
embedding $\sigma :M_{i}\hookrightarrow \mathcal{B}$ such that $\tau
^{-1}\circ \sigma $ is an identity on $K$ if and only if $\tau (f_{i})$ has
a root in $\mathcal{B}$. Indeed, suppose there exists $\gamma \in \mathcal{B}
$ such that $\gamma $ is a root of $\tau (f_{i})$. Then $\beta =\tau
^{-1}(\gamma )\in M_{X}$ and $f_{i}$ has a root $\beta \in M_{X}$. We claim
that
\begin{equation*}
K(\beta )=K(\alpha _{i})=M_{i}\text{.}
\end{equation*}%
Suppose otherwise. Then consider $\lambda :K(\alpha _{i})\rightarrow K(\beta
)$ keeping $K$ fixed, and its extension $\tilde{\lambda}$ to $M$. Since $%
\tilde{\lambda}$ fixes $K$, by our assumption on $\{M_{i}\}$, we conclude
that either $K(\beta )=K(\alpha _{i})$, or $K(\beta )\not\subset M$ and, in
particular, $K(\beta )\not\subset M_{X}$. Thus, $i\in X$ and $M_{i}\subset
M_{X}$, and we can set $\sigma (M_{i})=\tau (M_{i})$. Conversely, if the
required $\sigma $ exists, then $\tau ^{-1}\circ \sigma :M_{i}\rightarrow
M_{X}$, and thus $i\in X$ and $\tau (\alpha _{i})$ will satisfy $\tau
(f_{i}) $. Therefore, we just need to check systematically all elements of $%
\mathcal{B}$ until we find a root for $\tau (f_{i})$ or $\tau (f_{j})$.
\end{enumerate}
\end{proof}

\begin{example}
Let $K=\mathbb{Q}$, let $\{p_{i},i\in \omega \}$ be the listing of rational
primes. Let $f_{i}(T)=T^{2}-p_{i}$. It is clear that in this case the
sequence $\{M_{i},i\in \omega \}$, where $M_{i}=\mathbb{Q}(\sqrt{p_{i}})$,
is stable and totally linearly disjoint over $\mathbb{Q}$, and thus all the
requirements of Theorem \ref{thm:2} are satisfied.
\end{example}

\begin{example}
Let $K=\mathbb{Q}(x)$, where $x$ is not algebraic over $\mathbb{Q}$. Let $%
\{p_{i},i\in \omega \}$ be the listing of rational primes. Let $M_{i}=%
\mathbb{Q}(x,\sqrt[p_{i}]{x^{2}+1})$. Then the sequence $\{M_{i},i\in \omega
\}$ is linearly disjoint over $\mathbb{Q}(x)$ and is stable with respect to $%
\mathbb{Q}(x)$. Hence Theorem \ref{thm:2} applies again.
\end{example}

\begin{example}
Let $\mathbb{F}_{p}$ be a field of $p$ elements for some rational prime $p$.
Let $\{p_{i},i\in \omega \}$ be a listing of rational primes as before. Let $%
\alpha _{i}$ be of degree $p_{i}$ over $\mathbb{F}_{p}$. Let $M_{i}=\mathbb{F%
}_{p}(\alpha _{i})$. Then $\{M_{i},i\in \omega \}$ is a stable sequence
totally linearly disjoint over $\mathbb{F}_{p}$, and the Theorem \ref{thm:2}
applies.
\end{example}

\begin{example}
Let $\mathbb{G}_{p}$ be a finite field. Let $x$ be transcendental over $%
\mathbb{G}_{p}$. Let $\{p_{i},i\in \omega \}$ be a listing of rational
primes as before. Let $\alpha _{i}$ be of degree $p_{i}$ over $\mathbb{G}%
_{p} $. Let $M_{i}=\mathbb{G}_{p}(\alpha _{i},x)$. Then $\{M_{i},i\in \omega
\}$ is totally linearly disjoint over $\mathbb{G}_{p}(x)$ and is stable with
respect to $\mathbb{G}_{p}(x)$.
\end{example}

\begin{example}
Let $K=\mathbb{Q}(x)$ or $K=\mathbb{F}_{p}(x)$, where $x$ is not algebraic
over $\mathbb{Q}$ or $\mathbb{F}_{p}$, respectively. Let $M_{i}=\mathbb{Q}(%
\sqrt{x^{2}+i})$ or $M=\mathbb{F}_{p}(\sqrt{x^{2}+i})$. Then $\{M_{i},i\in
\omega \}$ is totally linearly disjoint over $\mathbb{Q}(x)$ and $\mathbb{F}%
_{p}(x)$, respectively, and is stable with respect to $\mathbb{Q}(x)$ and $%
\mathbb{F}_{p}(x)$, respectively.
\end{example}

\section{Rings whose isomorphism types have arbitrary Turing degrees}

\label{sec:rings} In this section we consider sequences of integrally closed
subrings of product formula fields: number fields (finite extensions of $%
\mathbb{Q}$), and finite extensions of rational function fields. In the case
of a function field, we will let the constant field be an arbitrary
computable field with a splitting algorithm. (For a discussion of fields
with splitting algorithms see Section 17.1, 17.2 of \cite{Fried}.)

Let $K$ be a product formula field. In the case of a function field, let $C$
be the constant field, and let $x\in K$ be a non-constant element. If
characteristic $p$ is such that $p>0$, then assume $x$ is not a $p$th power
in $K$. Under our assumptions, $K/\mathbb{Q}$ in the case of a number field,
or $K/C(x)$ in the case of a function field is a finite and separable
extension. Let $R=\mathbb{Z}$ in the case of a number field, and let $R=C[x]$
in the case of a function field. Then all prime ideals of $R$ correspond to
prime numbers or irreducible monic polynomials. In the case of $\mathbb{Z}$,
these ideals also represent all non-archimedean valuations of $\mathbb{Q}$,
while in the case of $C(x)$, a valuation corresponding to the degree of
polynomials does not correspond to a prime ideal of $C[x]$. (It does,
however, correspond to a prime ideal of $C[\frac{1}{x}]$.) Now, let $O_{K}$
be the integral closure of $\mathbb{Z}$ or $C[x]$ in $K$. Then $O_{K}$ is
called the ring of algebraic integers or integral functions, depending on
the choice of $K$. Now, the prime ideals of $R$ do not necessarily remain
prime in $O_{K}$, but every prime ideal of $R$ will have finitely many
factors in $O_{K}$. The set of all prime ideals of $O_{K}$ (together with
the factors of the degree valuation in the case of function fields) will
constitute what is called the set of \emph{primes of}\textit{\ }$K$. If $%
x\in (O_{K})^{\ast }$ and $\mathsf{p}$ is a prime ideal in $O_{K}$, then we
define $ord_{\mathsf{p}}z$ to be the largest nonnegative number such that $%
z\in (\mathsf{p})^{n}$. If $w\in K$, we write $w=z_{1}/z_{2}$ for some $%
z_{1},z_{2}\in (O_{K})^{\ast }$, and we let $ord_{\mathsf{p}}w=ord_{\mathsf{p%
}}z_{1}-ord_{\mathsf{p}}z_{2}$. Finally, we set $ord_{\mathsf{p}}0=\infty $.
For more material on valuations and primes of number fields and function
fields the reader is referred to \cite{C}, \cite{Januz}, \cite{Fried}, and
\cite{Riben}.

Using the Strong Approximation Theorem (see p. 21 of \cite{Fried} and p. 268
of \cite{Riben}), it can be shown that any integrally closed subring of $K$,
whose fraction field is $K$ is of the form

\begin{equation*}
O_{K,\mathcal{W}}=\{z\in K:(\forall \mathsf{p}\not\in \mathcal{W})[ord_{%
\mathsf{p}}z\geq 0]\}\text{,}
\end{equation*}%
where $\mathcal{W}$ is an arbitrary set of (non-archimedean) primes of $K$.
In the case $\mathcal{W}$ is finite, this ring is called a ring of $\mathcal{%
W}$-integers. Unfortunately, there is no universally accepted name for these
rings when $\mathcal{W}$ is infinite. Since $R$ is a principal ideal domain,
$O_{K}$ is a free $R$-module. The basis of $O_{K}$ as a free $R$-module is
called an integral basis of $K$ over $R$. Therefore, if we represent
elements of $K$ using their coordinates with respect to an integral basis, $%
O_{K}$ is computable under such representation.\bigskip

In what follows, we use an effective listing of primes of a product formula
field, and basically identify the $i$th prime on the list with a natural
number $i$ (for the purpose of satisfying the requirements of Theorem \ref%
{thm:1}). To carry out this plan, we need some way to represent the primes
and to make sure that the order at a given prime is computable. There are
several ways to do this. In the next lemma we describe one of them.

\begin{lemma}
\label{le:pr} Let $K$ be a product formula field. Consider the pairs of $K$%
-elements representing the primes of $K$, as described below. Then, given an
element $x\in K$ presented by its coordinates with respect to some integral
basis of $K$ over its rational subfield, there is an effective procedure to
determine which primes of $K$ occur in the divisor of $x$, as well as the
order of $x$ at each of the primes occurring in its divisor.
\end{lemma}

\begin{proof}
Let $\mathsf{p}$ be a prime of $K$. First, by the Weak Approximation
Theorem, there exists $t_{\mathsf{p}}\in O_{K}$ such that $ord_{\mathsf{p}%
}t_{\mathsf{p}}=1$ and for any $\mathsf{q}\not=\mathsf{p}$, conjugate to $%
\mathsf{p}$ over the rational field, $ord_{\mathsf{q}}t_{\mathsf{p}}=0$. We
will identify each prime of $K$ with a pair $(t_{\mathsf{p}},p)$, where $p$
is the prime number or the monic irreducible polynomial below $\mathsf{p}$
in the rational field.

Since we assume that the constant field is computable and has a splitting
algorithm, we can produce a computable sequence of all primes of the
rational field. Then for all but finitely many primes $p$ of the rational
field, using Lemma 4.1 on p. 135 of \cite{Sh}, Lemma 17.5 on p. 232 of \cite%
{Fried}, and Proposition 25 on p. 27 of \cite{L}, we have an effective
procedure for determining the number and relative degree of all $K$-factors
of $p$. Next, by a systematic search of $O_{K}$, using the rational norms of
elements of $O_{K}$, we can locate $t_{\mathsf{p}}$ for each factor of a
given rational prime $p$. Finally, for an arbitrary element $x$ of $K$, by
looking at its rational norm, we can determine a finite superset of the
primes occurring in its divisor. Then, using the pairs for these potentially
occurring primes, we can determine the actual divisor of $x$.
\end{proof}

\begin{remark}
For future use we also note here that for any $K$-prime $\mathsf{p}$, by the
Strong Approximation Theorem, there exists an element $z_{\mathsf{p}}\in K$
such that $\mathsf{p}$ is the only prime of $K$ with $ord_{\mathsf{p}}z_{%
\mathsf{p}}<0$. Furthermore, there exist $a_{\mathsf{p}},b_{\mathsf{p}}\in K$
such that $z_{\mathsf{p}}=\frac{a_{\mathsf{p}}}{b_{\mathsf{p}}}$. Having
constructed sequence $\{(p,t_{\mathsf{p}})\}$ for each prime of $K$, we can
effectively locate $z_{\mathsf{p}}$, $a_{\mathsf{p}}$, and $b_{\mathsf{p}}$.
\end{remark}

\medskip

We are now ready to state the main theorem of this section.

\begin{theorem}
\label{thm:4} Let $K$ be a computable product formula field with a finite
degree separable rational subfield $F$. $($So either $F=\mathbb{Q}$, or $%
F=C(x)$ for some $x\in K\setminus C$.$)$ In the case of a function field,
assume that the constant field $C$ is computable, finitely generated, and
has a splitting algorithm. Let $\{p_{i},i\in \omega \}$ be an effective
listing of primes of $F$ $($i.e., prime numbers or monic irreducible
polynomials in $x$ over $C)$. For each $p_{i}$, choose one $K$-factor $%
\mathsf{p}_{i}$, and let $\{(p_{i},t_{\mathsf{p}_{i}})\}$ be the listing of
primes of $K$ corresponding to $\{p_{i},i\in \omega \}$, where $t_{\mathsf{p}%
_{i}}$ has the least code in the set $\{t_{\mathsf{q}}:\mathsf{q}$ is a $K$%
-factor of $p\}$. Let $\mathcal{A}=R$ $($that is, $\mathcal{A}=\mathbb{Z}$
if $K$ is a number field, and $\mathcal{A}=C[x]$ if $K$ is a function field$%
) $. Let $\mathcal{W}=\emptyset $ if $K$ is a number field, and let $%
\mathcal{W}$ be the $($finite$)$ set of all $K$-poles of $x$ if $K$ is a
function field. For any $X\subseteq \omega $, let
\begin{equation*}
\mathcal{W}_{X}=\mathcal{W}\cup \{\mathsf{p}_{i}:i\in X\}\text{.}
\end{equation*}%
Let
\begin{equation*}
\mathcal{A}_{X}=O_{K,\mathcal{W}_{X}}\text{.}
\end{equation*}%
Then the sequence $\{\mathcal{A}_{i},i\in \omega \}$ satisfies the
conditions of Theorem \ref{thm:1}.
\end{theorem}

\begin{proof}
In the case of $K$ being a function field, we observe that for all $i\in
\omega {}$, we have that $\mathsf{p}_{i}$ is not a pole of $x$ in $K$. This
is true since each $\mathsf{p}_{i}$ is a factor of a $C(x)$-prime
corresponding to a polynomial in $x$. We next verify that all conditions of
Theorem \ref{thm:1} hold.

\begin{enumerate}
\item We show that $O_{K,\mathcal{W}_{X}}\leq _{T}X$ for any $X\subseteq
\omega $. Assuming we represent elements of $K$ as $n$-tuples of their
coordinates with respect to some integral basis of $K$ over its rational
subfield $F$, by Lemma \ref{le:pr}, we can effectively compute the divisors
of elements of $K$. Assuming the element is in $O_{K,\mathcal{W}_{\omega }}$
to begin with, we can next determine for each pair $(p_{i},\mathsf{p}_{i})$
with $p_{i}$ occurring in the norm, whether $i\in X$, and thus establish
whether the element is in $O_{K,\mathcal{W}_{X}}$. (Conversely, if we have
the characteristic function of $O_{K,\mathcal{W}_{X}}$, to determine whether
$i\in X$, it is enough to determine whether $z_{\mathsf{p}_{i}}\in O_{K,%
\mathcal{W}_{X}}$. Hence we have $O_{K,\mathcal{W}_{X}}\equiv _{T}X$.)

\item Let $\sigma :O_{K,\mathcal{W}{}_{i}}\hookrightarrow O_{K,{}\mathcal{W}%
_{X}}$ for some $i\in \omega $ and $X\subset \omega $ with $\sigma
|_{C[x]}=id$\ \ in the case $K$ is a function field. We show that $i\in X$.
First, observe that $\sigma $ can be extended to an automorphism of $K$
keeping $C(x)$ fixed. Next, observe that $\sigma (z_{\mathsf{p}_{i}})=\sigma
(a_{\mathsf{p}_{i}}/b_{\mathsf{p}_{i}})$ is an element that has negative
order only at a prime $\sigma (\mathsf{p}_{i})$. By the definition of $%
O_{K,{}\mathcal{W}_{X}}$, we have that $\sigma (\mathsf{p}_{i})\in \mathcal{W%
}_{X}\subset \mathcal{W}$. However, $\mathcal{W}_{X}$ contains only one
factor for each $p_{i}$. Therefore,
\begin{equation*}
\sigma (\mathsf{p}_{i})=\mathsf{p}_{i}\in \mathcal{W}_{X}\Leftrightarrow
i\in X\text{.}
\end{equation*}%
The case of $K$ being a number field is similar.

\item Let $\tau :O_{K,\mathcal{W}_{X}}\longleftrightarrow \mathcal{B}$ be an
isomorphism of rings. Let $i,j\in \omega $, $i\not=j$, be such that for
exactly one of $k=i$ and $k=j$, there exists an embedding $\sigma :O_{K,%
\mathcal{W}_{k}}\hookrightarrow \mathcal{B}$ such that $\tau ^{-1}\circ
\sigma |_{C[x]}=id$ in the case when $K$ is a function field. We show that,
relative to $\mathcal{B}$, we can compute which of the rings $O_{K,\mathcal{W%
}_{i}}$ and $O_{K,\mathcal{W}_{j}}$ is embeddable into $\mathcal{B}$ in the
prescribed manner. Since $C$ is finitely generated and $O_{K}$ has a basis
over $R$, in $\mathcal{B}$ we can list $\tau (O_{K})$ and effectively find $%
\tau (a_{\mathsf{p}_{i}})$, $\tau (b_{\mathsf{p}_{i}})$, $\tau (a_{\mathsf{p}%
_{j}})$, and $\tau (b_{\mathsf{p}_{j}})$. Finally, we systematically look
for the solution in $\mathcal{B}$ of the equations $\tau (b_{\mathsf{p}%
_{i}})Z=\tau (a_{\mathsf{p}_{i}})$, and $\tau (b_{\mathsf{p}_{j}})W=\tau (a_{%
\mathsf{p}_{j}})$. It is clear that, by assumption, exactly one of these
equations will have a solution. On the other hand, $\tau (b_{\mathsf{p}%
_{k}})Z=\tau (a_{\mathsf{p}_{k}})$ has a solution in $\mathcal{B}$ if and
only if $\sigma :O_{K,{}\mathcal{W}_{k}}\hookrightarrow \mathcal{B}$ as
specified above exists.
\end{enumerate}
\end{proof}

In lieu of examples illustrating this theorem, which are pretty easy to
generate, we offer a different version of Theorem \ref{thm:4}. The proof of
this version is very similar to the proof of Theorem \ref{thm:4}.

\begin{theorem}
\label{thm:4.1} Let $K$ be a computable product formula field with a finite
degree separable rational subfield $F$. In the case of a function field,
assume that the constant field $C$ is computable, finitely generated, and
has a splitting algorithm. Let $\{\mathsf{p}_{i},i\in \omega \}$ be an
effective listing of primes of $K$, excluding in the function field case the
poles of some non-constant element $x$ such that $F=C(x)$. Let $\mathcal{A}%
=O_{K}$. Let $\mathcal{W}=\emptyset $ if $K$ is a number field, and let $%
\mathcal{W}$ be the $($finite$)$ set of all $K$-poles of $x$. For any $%
X\subseteq \omega $, let
\begin{equation*}
\mathcal{W}_{X}=\bigcup_{i\in X}\{\mathsf{p}\text{
is a factor of }p_{i}\} \cup \mathcal{W} \text{.}
\end{equation*}%
Let $\mathcal{W}_{i}=\mathcal{W}\cup \{$all $K$-factors of $p_{i}\}$. Let
\begin{equation*}
\mathcal{A}_{X}=O_{K,\mathcal{W}_{X}}.
\end{equation*}%
Let $\mathcal{A}_{i}=O_{K,\mathcal{W}_{i}}$. Then the sequence $\{\mathcal{A}%
_{i},i\in \omega \}$ satisfies the conditions of Theorem \ref{thm:1}.
\end{theorem}

\section{Torsion-free abelian groups of arbitrary finite rank and of
arbitrary Turing degree}

\label{sec:groups}Rank $1$, torsion-free, abelian groups are isomorphic to
subgroups of $(\mathbb{Q},+)$. There is a known classification, due to Baer
\cite{baer}, of countable, rank $1$, torsion-free, abelian groups. The
account here will generally follow the one in a book  by Fuchs \cite%
{Fuchs}, Volume $2$. Given such a group $\mathcal{G}$, for any prime $p$, we define a
function $h_{p}:G\rightarrow \mathbb{\omega }$ by setting $h_{p}(a)$ equal
to the largest natural number $k$ such that there is some $b\in G$ with $%
p^{k}b=a$. If no such $k$ exists, we set $h_{p}(a)=\infty $. Now define the
\emph{characteristic} of $a$ to be the sequence
\begin{equation*}
\chi _{\mathcal{G}}(a)=(h_{p_{1}}(a),h_{p_{2}}(a),\dots )\text{,}
\end{equation*}%
where $(p_{i})_{i\in \omega }$ is a list of all prime numbers.

In some torsion-free abelian groups (for example, $(\mathbb{Q},+)$), it is
the case that all nonzero elements have the same characteristic. In these
groups, we would need to look no further for invariants. However, in some
others (for example, $(\mathbb{Z},+)$) the characteristics of the various
elements are essentially the same, but not identical. We say that two
characteristics are equivalent if they are equal except in a finite number
of places, and in all places where they differ, both are finite. An
equivalence class of characteristics under this relation is called a \emph{%
type}. If $\chi _{\mathcal{G}}(a)$ belongs to a type $\mathbf{t}$, then we
set $\mathbf{t}_{\mathcal{G}}(a)=\mathbf{t}$, and say that $a$ is of type $%
\mathbf{t}$. A group $\mathcal{G}$ in which any two nonzero elements have
the same type $\mathbf{t}$ is \emph{homogeneous}, and we say that $\mathbf{t}%
(\mathcal{G})=\mathbf{t}$ is the type of $\mathcal{G}$. In particular, we
note that any torsion-free abelian group of rank $1$ is homogeneous.

\begin{proposition}
\label{Baer}$(${Baer \cite{baer}}${)}$ If $\mathcal{G}$ and $\mathcal{H}$
are torsion-free abelian groups of rank $1$, then $\mathcal{G}\simeq
\mathcal{H}$ if and only if $\mathbf{t}(\mathcal{G})=\mathbf{t}(\mathcal{H})$%
.
\end{proposition}

Knight and Downey established that for an arbitrary Turing degree $\mathbf{d}
$, there exists a torsion-free abelian group of rank $1$ and finite type,
whose degree of the isomorphism type is $\mathbf{d}$ (see \cite{D}). Downey
and Jockusch showed in \cite{D} that even rank $1$ groups of finite type can
fail to have degree. On the other hand, it follows from a
computability-theoretic result of Coles, Downey, and Slaman, using a result
by Downey, Jockusch, and Solomon, that every torsion-free abelian group of
rank $1$ has first jump degree (see \cite{CDS}). Coles, Downey, and Slaman
proved that for every set $A\subseteq \omega $, there is a Turing degree
that is the least degree of the jumps of all sets $X$ for which $A$ is
computably enumerable in $X$. On the other hand, Richter \cite{R}\
constructed a non-computably enumerable set that is computably enumerable in
two sets that form a minimal pair. This construction implies that there is a
set $A$ such that the set of all $X$ for which $A$ is computably enumerable
in $X$ has no member of the least Turing degree.

Knight-Downey's proof of the claim that for an arbitrary Turing degree $%
\mathbf{d}$, there exists a torsion-free abelian group of rank $1$ (and
finite type) whose degree of the isomorphism type is $\mathbf{d}$ is as
follows. Let a set $D\subseteq \omega $ be of degree $\mathbf{d}$. Let $%
\mathcal{G}$ be a rank $1$ group defined by the type sequence $(a_{n})_{n\in
\omega }$ such that
\begin{equation*}
a_{n}=\left\{
\begin{array}{cc}
1 & \text{if }n\in D\oplus \bar{D}\text{,} \\
0 & \text{if }n\notin D\oplus \bar{D}\text{.}%
\end{array}%
\right.
\end{equation*}%
There is an isomorphic copy of $\mathcal{G}$ computable in $D$. Conversely,
let $\mathcal{H}\cong \mathcal{G}$. Then $\mathcal{H}$ has the type sequence
$(a_{n})_{n\in \omega }$ (by Proposition \ref{Baer}). The type sequence of $%
\mathcal{H}$ is computably enumerable in $\mathcal{H}$, and hence $D\oplus
\bar{D}$ is computably enumerable in $\mathcal{H}$. Thus $D\leq _{T}\mathcal{%
H}$.

On the other hand, the construction from Section \ref{sec:rings} can be
adapted to give examples of torsion-free abelian groups of finite rank
satisfying conditions of Theorem \ref{thm:1}, and thus produce torsion-free
abelian groups of arbitrary rank and of arbitrary Turing degree $\mathbf{d}$.

\begin{theorem}
\label{thm:groups} Let $k$ be a positive integer. Let $\mathcal{C}$ be the
class of all subgroups of $\mathbb{Q}^{k}$. Let $\{p_{i},i\in \omega \}$ be
a listing of rational primes. Let $\mathcal{A}=\mathbb{Z}^{k}$. For any $%
X\subseteq \omega $, let $\mathcal{V}_{X}=\{p_{i}:i\in X\}$. Let
\begin{equation*}
\mathcal{A}_{X}=(O_{\mathbb{Q},_{\mathcal{V}_{X}}})^{k}.
\end{equation*}
For each $i\in \omega $, let $\mathcal{A}_{i}=(O_{\mathbb{Q},\{p_{i}\}})^{k}$%
. Then the sequence $\{\mathcal{A}_{i},i\in \omega \}$ satisfies the
conditions of Theorem \ref{thm:1}.
\end{theorem}

The proof of this theorem is almost identical to the proof of Theorem \ref%
{thm:4}.

\section{Algebraic structures whose isomorphism types have no Turing degrees}

We now want to investigate structures whose isomorphism types have no Turing
degrees. The enumeration reducibility will play an important role. First, we
define the canonical index $u$ of a finite set $D_{u}\subset \omega $. Let $%
D_{0}=\emptyset $. For $u>0$, let $D_{u}=\{n_{0},\ldots ,n_{k-1}\}$, where $%
n_{0}<\ldots <n_{k-1}$ and $u=2^{n_{0}}+\ldots +2^{n_{k-1}}$. For sets $%
X,Y\subseteq \omega $, we say that $Y$ is \emph{enumeration reducible} to $X$%
, in symbols $X\leq _{e}Y$, if there is a computably enumerable binary
relation $E$ such that

\begin{equation*}
(x\in X)\Leftrightarrow (\exists u)[D_{u}\subseteq Y\wedge E(x,u)]\text{.}
\end{equation*}%
Equivalently, $X\leq _{e}Y$ iff for every set $S$, if $Y$ is computably
enumerable relative to $S$, then $X$ is computably enumerable relative to $S$%
. Richter showed that a modification of Theorem 2.1 on p. 725 in \cite{R},
obtained by replacing Turing reducibility in $\mathcal{A}_{X}\leq _{T}X$ by
the enumeration reducibility $\mathcal{A}_{X}\leq _{e}X$, yields a very
different conclusion, Theorem 2.3 on p. 726 in \cite{R} --- that there is a
set $X$ such that the isomorphism type of $\mathcal{A}_{X}$ does not have a
degree. Richter used this theorem to show that there is a torsion abelian
group whose isomorphism type does not have a degree. The following result is
a generalization of Richter's Theorem 2.3. It should be noted that the
hypotheses of this theorem are similar to those of Theorem \ref{thm:1}.

\begin{theorem}
Let $\mathcal{C}$ be a class of countable structures in a finite language $L$%
, closed under isomorphisms. Assume that there is a computable sequence $\{%
\mathcal{A}_{i},i\in \omega \}$ of computable $($possibly infinite$)$
structures in $\mathcal{C}$ satisfying the following conditions.

\begin{itemize}
\item There exists a finitely generated structure $\mathcal{A}\in \mathcal{C}
$ such that for all $i\in \omega $, we have that $\mathcal{A}\subset
\mathcal{A}_{i}$.

\item For any $X\subseteq \omega $, there is a structure $\mathcal{A}_{X}$
in $\mathcal{C}$ such that $\mathcal{A}\subset \mathcal{A}_{X}$ and
\begin{equation*}
\mathcal{A}_{X}\leq _{e}X\text{,}
\end{equation*}%
and for every $i\in \omega $, there exists an embedding $\sigma $ such that
\begin{equation*}
\sigma :\mathcal{A}_{i}\hookrightarrow \mathcal{A}_{X}\text{, \ }\sigma _{|%
\mathcal{A}}=id\text{,}
\end{equation*}%
iff $i\in X$.

\item Suppose a structure $\mathcal{B}$ is such that for some $X\subseteq
\omega $ we have that $\mathcal{B\simeq A}_{X}$ under isomorphism $\tau :%
\mathcal{A}_{X}\longleftrightarrow \mathcal{B}$. Then from any enumeration
of $\mathcal{B}$, we can effectively pass to an enumeration of the set of
indices $i$ such that $\mathcal{A}_{i}\hookrightarrow \mathcal{B}$ under an
embedding $\sigma $ such that $(\tau ^{-1}\circ \sigma )_{|\mathcal{A}}=id$.
\end{itemize}

Then there is a structure in $\mathcal{C}$ whose isomorphism type has no
Turing degree.
\end{theorem}

\begin{proof}
The proof follows from the following well-known lemma.

\begin{lemma}
\label{setwnoleast} There is a set $X\subseteq \omega $ such that the set of
functions
\begin{equation*}
\{f:ran(f)=X\}
\end{equation*}%
has no least Turing degree.
\end{lemma}

We omit the proof of Lemma \ref{setwnoleast} here. The idea is to construct
a set $X$, which is not computably enumerable, but the set of enumerations
of $X$ contains two enumerations whose Turing degrees form a minimal pair. A
full proof may be found in \cite{R}.

Let $X$ be as in Lemma \ref{setwnoleast}, and, toward contradiction, let $%
\mathcal{M}$ be a copy of $\mathcal{A}_{X}$ with least Turing degree. Now,
we will argue that $\mathcal{M}$ has a least enumeration. Let $\nu $ be an
enumeration of $\mathcal{M}$ where $deg(\nu )=deg(\mathcal{M})$. Let $g$ be
another enumeration of $\mathcal{M}$ with $\nu \nleq _{T}g$. By padding
(i.e.\ by passing to a structure whose elements are of the form $(x,t)$,
where $x$ is enumerated into $M$ at time $t$, and where the operations are
the obvious ones), we can obtain an isomorphic copy of $\mathcal{M}$ with
the same Turing degree as $g$, which is a contradiction. Hence, by the final
assumption of the theorem, we can pass effectively from $\nu $ to an
enumeration $\nu _{\mathcal{M}}$ of $X$.

Let $f$ be an enumeration of $X$, and we will show that $\nu _{\mathcal{M}%
}\leq _{T}f$. By assumption, we can pass effectively from $f$ to an
enumeration $g$ of $\mathcal{A}_{X}$, and, since $\mathcal{A}_{X}\leq _{e}X$
and $\mathcal{M}\leq _{T}\mathcal{A}_{X}$ (because here we identify $%
\mathcal{M}$ with its atomic diagram), we can pass effectively from $g$ to
an enumeration $\tilde{g}$ of $\mathcal{M}$. Now $\mathcal{M}\leq _{T}\tilde{%
g}$, so we have $\nu _{\mathcal{M}}\leq _{T}\mathcal{M}\leq _{T}f$, as was
to be shown. This is a contradiction.
\end{proof}

\begin{theorem}
\label{thm:enum} There are countable fields, rings, and torsion-free abelian
groups of arbitrary finite rank, whose isomorphism types do not have Turing
degrees.
\end{theorem}

\begin{proof}
Since the statements of Theorem \ref{thm:enum} and Theorem \ref{thm:1}
differ in the place where $\mathcal{A}_{X}\leq _{T}X$ is replaced with $%
\mathcal{A}_{X}\leq _{e}X$, to prove the analogs of Theorems \ref{thm:2}, %
\ref{thm:4}, and \ref{thm:groups}, we need to show that the respective
structures satisfy $\mathcal{A}_{X}\leq _{e}X$.

We will start with the field case, where all the notation and assumptions
are as in Theorem \ref{thm:2}. We need to show that $M_{X}\leq _{e}X$. Let $%
\phi :\mathbb{\omega }\longrightarrow X$ be any listing of $X$. Then, using $%
\phi $, we can list the sets $\{\alpha _{\phi (i)}\}$ and $\{\alpha _{\phi
(i_{1})}^{m(\phi (i_{1}))}\ldots \alpha _{\phi (i_{k})}^{m(\phi (i_{k}))}\}$%
, where $0\leq m(\phi (i_{l}))<\deg (f_{\phi (i_{l})})$ and $k\in \mathbb{%
\omega }$. Thus we will be able to list the basis of $M_{X}$ over $K$, and
then $M_{X}$ itself. (As in the case of Turing reducibility, it is not hard
to see that $X\leq _{e}M_{X}$ also, and therefore, we really have
enumeration equivalence.)

We now proceed to the ring case. Here, all the notation and assumptions are
as in Theorem \ref{thm:4}. Let $\phi $ be as above, and note that, given $%
\phi $, we can list the set $\{p_{\phi (i)},\mathsf{p}_{\phi (i)}\}$. The
listing of primes will then allow the listing of $O_{K,\mathcal{W}_{X}}$,
where we would proceed by testing elements of $K$ to see if the primes in
the denominator of their divisors have appeared in the list already. This
testing process is effective as discussed in Section \ref{sec:rings}. (As
above, we also have here that $X\leq O_{K,{}\mathcal{W}_{X}}$.) The case of
the abelian groups from Section \ref{sec:groups} is again almost identical
to the case of the rings.

For the condition on listing the $i$ such that $\mathcal{A}_i$ is embeddable
in a copy of $\mathcal{A}_X$, the earlier proofs suffice. In each case, we
proved the hypotheses for Theorem \ref{thm:1} by first establishing the
hypotheses for this theorem.
\end{proof}


\bigskip

\begin{thebibliography}{99}
\bibitem{A} \textsc{E. Artin}, \textbf{Algebraic Numbers and Algebraic
Functions}, Gordon and Breach Science Publishers, New York, 1967.

\bibitem{AJK} \textsc{C. J. Ash,\ C. G. Jockusch, Jr., }and\textsc{\ J. F.
Knight}, \textit{Jumps of orderings}, \textbf{Transactions of the American
Mathematical Society} 319 (1990), pp. 573--599.

\bibitem{AK} \textsc{C. J. Ash }and\textsc{\ J. F. Knight}, \textbf{%
Computable Structures and the Hyperarithmetical Hierarchy}\textit{,}
Elsevier, Amsterdam, 2000.

\bibitem{baer} \textsc{R. Baer}, \textit{Abelian groups without elements of
finite order}, \textbf{Duke Mathematical Journal} 3 (1937), pp. 68--122.

\bibitem{C} \textsc{C. Chevalley}, \textbf{Introduction to the Theory of
Algebraic Functions of One Variable}, Mathematical Surveys VI, American
Mathematical Society, New York, 1951.

\bibitem{Csima} \textsc{B. Csima}, \textit{Degree spectra of prime models},
\textbf{Journal of Symbolic Logic} 69 (2004), pp. 430--442.

\bibitem{CDS} \textsc{R. J. Coles, R. G. Downey, }and\textsc{\ T. A. Slaman,}
\textit{Every set has a least jump enumeration}, \textbf{Journal of the
London Mathematical Society} 62 (2000), pp. 641--649.

\bibitem{D} \textsc{R. G. Downey}, \textit{On presentations of algebraic
structures}, \textbf{Complexity, Logic, and Recursion Theory} (A. Sorbi,
editor), Lecture Notes in Pure and Applied Mathematics 187, Marcel Dekker,
New York, 1997, pp. 157--205.

\bibitem{DK} \textsc{R. Downey }and\textsc{\ J. F. Knight}, \textit{\
Orderings with }$\alpha $\textit{th jump degree }\textbf{0}$^{(\alpha )}$,
\textbf{Proceedings of the American Mathematical Society} 114 (1992), pp.
545--552.

\bibitem{Fried} \textsc{M. D. Fried }and\textsc{\ M. Jarden}, \textbf{Field
Arithmetic}, first edition, Springer-Verlag, Berlin, 1986.

\bibitem{Fuchs} \textsc{L. Fuchs}, \textbf{Infinite Abelian Groups},\
Academic Press, New York, 1973.

\bibitem{H1} \textsc{V. S. Harizanov}, \textit{Computability-theoretic
complexity of countable structures}, \textbf{Bulletin of Symbolic Logic} 8
(2002), pp. 457--477.

\bibitem{H} \textsc{V. S. Harizanov}, \textit{Pure computable model theory, }
\textbf{Handbook of Recursive Mathematics}, vol. 1 (Yu. L. Ershov, S. S.
Goncharov, A. Nerode and J. B. Remmel, editors, V. W. Marek, assoc. editor),
Elsevier, Amsterdam, 1998, pp. 3--114.

\bibitem{HM} \textsc{V. Harizanov }and\textsc{\ R. Miller}, Spectra of
structures and relations, submitted.

\bibitem{HKSS} \textsc{D. R. Hirschfeldt, B. Khoussainov, R. A. Shore, }and
\textsc{\ A. M. Slinko}, \textit{Degree spectra and computable dimension in
algebraic structures}, \textbf{Annals of Pure and Applied Logic} 115 (2002),
pp. 71--113.

\bibitem{Januz} \textsc{G. J. Janusz}, \textbf{Algebraic Number Fields},
Pure and Applied Mathematics 55, Academic Press, New York, 1973.

\bibitem{JS1} \textsc{C. G. Jockusch, Jr. }and\textsc{\ R. I. Soare},
\textit{\ Boolean algebras, Stone spaces, and the iterated Turing jump},
\textbf{Journal of Symbolic Logic} 59 (1994), pp. 1121--1138.

\bibitem{Kh} \textsc{A. N. Khisamiev}, \textit{On the upper semilattice} $%
\mathcal{L}_{E}$, \textbf{Siberian Mathematical Journal} 45 (2004), pp.
211--228 (in Russian); pp. 173--187 (English translation).

\bibitem{K} \textsc{J. F. Knight}, \textit{Degrees coded in jumps of
orderings}, \textbf{Journal of Symbolic Logic} 51 (1986), pp. 1034--1042.

\bibitem{L} \textsc{S. Lang}, \textbf{Algebraic Number Theory},
Addison-Wesley, Reading, 1970.

\bibitem{Marker} \textsc{D. Marker}, \textbf{Model Theory}, An Introduction,
Springer-Verlag, New York, 2002.

\bibitem{Millar} \textsc{T. S. Millar}, \textit{Pure recursive model theory}%
, \textbf{Handbook of Computability Theory} (E. R. Griffor, editor),
Elsevier, Amsterdam, 1999, pp. 507--532.

\bibitem{O} \textsc{S. Oates}, \textbf{Jump Degrees of Groups}, PhD
dissertation, University of Notre Dame, 1989.

\bibitem{Rabin} \textsc{M. O. Rabin}, \textit{Computable algebra, general
theory and theory of computable fields}, \textbf{Transactions of the
American Mathematical Society} 95 (1960), pp. 341--360.

\bibitem{Riben} \textsc{P. Ribenboim}, \textbf{The theory of classical
valuations}, Springer Monographs in Mathematics, Springer-Verlag, New York,
1999.

\bibitem{R} \textsc{L. J. Richter}, \textit{Degrees of structures}, \textbf{%
Journal of Symbolic Logic} 46 (1981), pp. 723--731.

\bibitem{Rogers} \textsc{H. Rogers, Jr.}, \textbf{Theory of Recursive
Functions and Effective Computability}, McGraw-Hill, New York, 1967.

\bibitem{Sh} \textsc{A. Shlapentokh}, \textit{Hilbert's Tenth Problem over
number fields, a survey}, \textbf{Hilbert's Tenth Problem: Relations with
Arithmetic and Algebraic Geometry}, (J. Denef, L. Lipshitz, T. Pheidas, and
J. Van Geel, editors), Contemporary Mathematics 270, American Mathematical
Society, Providence, 2000, pp. 107--137.

\bibitem{S} \textsc{R. I. Soare}, \textbf{Recursively Enumerable Sets and
Degrees,} Springer-Verlag, Berlin, 1987.
\end{thebibliography}
\end{document}